\documentclass[12pt,a4paper,reqno]{amsart}
\usepackage{amsfonts}
\usepackage{graphicx}   
\usepackage{graphics}
\usepackage{epsfig}
\usepackage{anysize}
\marginsize{2.5cm}{2.5cm}{2.5cm}{2.5cm}

\def\({\left(}
\def\){\right)}

\begin{document}
\setcounter{page}{1}

\title []
{Investigation of solutions of boundary value problems for a
composite type equation with non-local boundary conditions}
\author{Aliev N.A., Aliev A.M.  }
\maketitle
\begin{center} \emph{Faculty of Applied Mathematics and Cybernetics,Institute of
Applied Mathematics,\\ Baku State University,Z.Khalilov
str.23,AZ1148, Azerbaijan,\\e-mail:nihan\_aliev@rambler.ru,
aahmad07@rambler.ru }\end{center}

\begin{abstract}
Since the order of elliptic type model equation (Laplace equation)
is two [1], [2], then it is natural the order of composite type
model equation must be [3] [4] [5] three. At each point of the
domain under consideration these equations have both real and
complex characteristics.

Notice that a boundary value problem for a composite type equation
of second order first appeared in the  paper [6].

The method for investigating the Fredholm property of boundary
value problems is distinctive and belongs to one of  the authors
of  the present paper.


\end{abstract}

\section*{Introduction.}

The paper is devoted to the investigation of boundary value problems for a composite type equation of second order.

This was possible as earlier we investigated an elliptic type boundary value problem of first order ( Cauchy-Riemann equation) for which the boundary is a carrier of boundary conditions [7]. Notice that in this case it is impossible to determine local boundary conditions  (since undeterminacy  is obtained). Therefore non-local boundary conditions were considered.

Necessary conditions that contain singular integrals are obtained proceeding from fundamental solution of Cauchy-Riemann equation [8]. Considering that we are on a spectrum, regularization of these singularities  is also conducted in distinctive way [6]. Joining regularized necessary conditions with the given boundary conditions we get a sufficient condition for Fredholm property of the stated boundary value problems.

Notice that in [8] the cited investigation of the process in a nuclear reactor leads to a boundary value problem for first order integro-differential equation in three-dimensional space where not all the space is a carrier of the given local boundary condition.

\section*{Problem statement }

 Let's consider the following boundary value problem:

\begin{equation} \label{GrindEQ__1_}
\ell \, u\equiv \frac{\partial ^{2} u(x)}{\partial x_{2}^{2} } +i\frac{\partial ^{2} u(x)}{\partial x_{1} \, \partial x_{2} } =0,\, \, \, \, \, x\in D,
\end{equation}

\begin{equation} \label{GrindEQ__2_}
\ell _{k} u\equiv \left. \frac{\partial u(x)}{\partial \, x_{2} } \right|_{x_{2} =\gamma _{k} (x_{1} )} +\alpha _{k} u(x_{1} ,\gamma _{k} (x_{1} ))=\varphi _{k} (x_{1} ),\, \, \, k=1,\, 2;\, \, x_{1} \in \left[a_{1} ,b_{1} \right]
\end{equation}
where $i=\sqrt{-1} $, $D\subset R^{2} $- is a bounded domain convex in the direction $x_{2} $, the boundary $\Gamma =\bar{D}\backslash D$ - is a Liapunov line, $\gamma _{k} (x_{1} )$, $k=1,\, 2$ are the equations of open lines $\Gamma _{k} $ ($\Gamma _{1} \bigcup \Gamma _{2} =\Gamma $), obtained from the boundary $\Gamma $ of the domain $D$ by means of orthogonal projection of this domain on the axis $x_{1} $ parallel to the axis  $x_{2} $ and  $\left[a_{1} ,\, b_{1} \right]=np_{x_{1} } \Gamma _{1} =np_{x_{1} } \Gamma _{2} $. In the given boundary conditions \eqref{GrindEQ__2_}  $\alpha _{k} $ ($k=1,\, 2$) are constants,  $\varphi _{k} (x_{1} )$, $k=1,\, 2$;  $x_{1} \in \left[a_{1} ,\, b_{1} \right]$  are sufficiently smooth functions. Boundary conditions \eqref{GrindEQ__2_} are assumed to be linear independent .

\section*{Fundamental solution.}

Applying the Fourier transform to equation \eqref{GrindEQ__1_} we get a fundamental solution in the form

\begin{equation} \label{GrindEQ__3_}
U(x-\xi )=\frac{-1}{4\pi ^{2} } \int _{R^{2} }\frac{e^{i(\alpha ,\, x-\xi )} }{\alpha _{2} (\alpha _{2} +i\, \alpha _{1} )} d\alpha  ,
\end{equation}
where

\[(\alpha _{1} x-\xi )=\sum _{j=1}^{2}\alpha _{j} (x_{j} -\xi _{j} ) .\]
Further, since

\[\frac{1}{2\pi i} \int _{R}\frac{e^{i\alpha _{2} (x_{2} -\xi _{2} )} }{\alpha _{2} } d\alpha _{2}  =e(x_{2} -\xi _{2} ),\]
where   $e(x_{2} -\xi _{2} )$  is a unique  symmetric Heaviside function, from  \eqref{GrindEQ__3_} we get:

\begin{equation} \label{GrindEQ__4_}
\frac{\partial U(x-\xi )}{\partial x_{2} } +i\frac{\partial U(x-\xi )}{\partial x_{2} } =e(x_{2} -\xi _{2} )\, \delta (x_{1} -\xi _{1} ).
\end{equation}
    Finally considering that \eqref{GrindEQ__3_} is a fundamental solution of the composite type equation \eqref{GrindEQ__1_}, we get that  $\frac{\partial U(x-\xi )}{\partial x_{2} } $  is a fundamental solution of the Cauchy-Riemann equation. Making negligible changes in the fundamental solution of the Cauchy-Riemann equation [8] we get a fundamental solution in the direction $x_{2} $

\begin{equation} \label{GrindEQ__5_}
\frac{\partial U(x-\xi )}{\partial x_{2} } =\frac{1}{2\pi } \frac{\theta (x_{2} -\xi _{2} )+\theta (\xi _{2} -x_{2} )}{x_{2} -\xi _{2} +i(x_{1} -\xi _{1} )} .
\end{equation}
here  $\theta (x_{2} -\xi _{2} )+\theta (x_{2} -\xi _{2} )=1$, if none differentiation operation is produced on it, since each addend has a break and contribution of this break appears in differentiation. Thus, for fundamental solution \eqref{GrindEQ__3_} of composite type equation \eqref{GrindEQ__1_} we get:

\begin{equation} \label{GrindEQ__6_}
U(x-\xi )=\frac{1}{2\pi } \int _{0}^{x_{2} }\frac{\theta (t-\xi _{2} )+\theta (\xi _{2} -t)}{t-\xi _{2} +i(x_{1} -\xi _{1} )} dt ,
\end{equation}
i.e. it holds the following statement:

\textbf{Theorem 1.} For a composite type equation of second order
\eqref{GrindEQ__1_} a fundamental solution in the direction $x_{2}
$ is of the form  \eqref{GrindEQ__6_}.

 This means that if we differentiate $U(x-\xi )$  twice with respect to $x_{2} $ and twice with respect to the mixed derivatives  $x_{1} $ and $x_{2} $, the Dirac delta function (two-dimensional) appears only in the derivative of second order with respect to $x_{2} $.

\section*{Necessary conditions.}

Multiplying equation  \eqref{GrindEQ__1_} by fundamental solution \eqref{GrindEQ__6_} and integrating if in domain $D$, applying Ostrogradskii-Gauss formula [8], we get formula similar to Green's second formula that after application of  fundamentality properties of function \eqref{GrindEQ__6_} get the form:

\[\int _{a_{1} }^{b_{1} }\left. \left[u(x)\frac{\partial U(x-\xi )}{\partial x_{2} } -\frac{\partial u(x)}{\partial x_{2} } U(x-\xi )\right]\right|_{x_{2} =\gamma _{1} (x_{1} )}^{\gamma _{2} (x_{1} )} dx_{1}  +i\int _{a_{1} }^{b_{1} }\left. u(x)\frac{\partial U(x-\xi )}{\partial x_{1} } \right|_{x_{2} =\gamma _{1} (x_{1} )}^{\gamma _{2} (x_{1} )} dx_{1}  +\]

\[+i\int _{a_{1} }^{b_{1} }\left. \frac{\partial u(x)}{\partial x_{2} } U(x-\xi )\right|_{x_{2} =\gamma _{2} (x_{1} )} \gamma '_{2} (x_{1} )dx_{1}  -i\int _{a_{1} }^{b_{1} }\left. \frac{\partial U(x)}{\partial x_{2} } U(x-\xi )\right|_{x_{2} =\gamma _{1} (x_{1} )} \gamma '_{2} (x_{1} )dx_{1}  =\]

\begin{equation} \label{GrindEQ__7_}
=\left\{\begin{array}{l} {u(\xi ),\, \, \, \, \, \, \, \xi \in D,} \\ {\frac{1}{2} u(\xi ),\, \, \xi \in \Im } \end{array}\right.
\end{equation}
The second expression in formula \eqref{GrindEQ__7_} is one of the necessary conditions. This condition has the form:

\[u(\xi _{1} ,\gamma _{1} (\xi _{1} ))=u(\xi _{1} ,\gamma _{2} (\xi
_{1} ))-2\int _{a_{1} }^{b_{1} }\left. \frac{\partial
u(x)}{\partial x_{2} } \right|_{x_{2} =\gamma _{2} (x_{1} )}
U(x_{1} -\xi _{1} ,\gamma _{2} (x_{1} )-\gamma _{1}
 (\xi _{1} ))\left[1-i\gamma '_{2} (x_{1} )\right]dx_{1}  +\]\\\begin{equation} \label{GrindEQ__8_}+2\int _{a_{1} }^{b_{1} }\left. \frac{\partial u(x)}{\partial x_{2} } \right|_{x_{2} =\gamma _{1} (x_{1} )} U(x_{1} -\xi _{1} ,\gamma _{1} (x_{1} )-\gamma _{1} (\xi _{1} ))\left[1-i\gamma '_{1} (x_{1} )\right]dx_{1}
\end{equation}
In exactly the same way to  [9]  and   [10],  we get the following necessary conditions:

\[\left. \frac{\partial u(\xi )}{\partial \xi _{1} } \right|_{\xi _{2} =\gamma _{1} (\xi _{1} )} =\left. \frac{\partial u(\xi )}{\partial \xi _{1} } \right|_{\xi _{2} =\gamma _{2} (\xi _{1} )} -i\left. \frac{\partial u(\xi )}{\partial \xi _{2} } \right|_{\xi _{2} =\gamma _{2} (\xi _{1} )} -\]

\[-2i\int _{a_{1} }^{b_{1} }\left. \frac{\partial u(x)}{\partial x_{2} } \right|_{x_{2} =\gamma _{1} (x_{1} )} \left. \frac{\partial U(x-\xi )}{\partial x_{2} } \right|_{\begin{array}{l} {x_{2} =\gamma _{1} (x_{1} )} \\ {\xi _{2} =\gamma _{1} (\xi _{1} )} \end{array}} \left[1-i\gamma '_{1} (x_{1} )\right]dx_{1}  +\]

\begin{equation} \label{GrindEQ__9_}
+2i\int _{a_{1} }^{b_{1} }\left. \frac{\partial u(x)}{\partial x_{2} } \right|_{x_{2} =\gamma _{2} (x_{1} )} \left. \frac{\partial U(x-\xi )}{\partial x_{2} } \right|_{\begin{array}{l} {x_{2} =\gamma _{2} (x_{1} )} \\ {\xi _{2} =\gamma _{1} (\xi _{1} )} \end{array}} \left[1-i\gamma '_{2} (x_{1} )\right]dx_{1}  ,
\end{equation}

\[\left. \frac{\partial u(\xi )}{\partial \xi _{2} } \right|_{\xi _{2} =\gamma _{1} (\xi _{1} )} =2\int _{a_{1} }^{b_{1} }\left. \frac{\partial u(x)}{\partial x_{2} } \right|_{x_{2} =\gamma _{2} (x_{1} )} \left. \frac{\partial U(x-\xi )}{\partial x_{2} } \right|_{\begin{array}{l} {x_{2} =\gamma _{2} (x_{1} )} \\ {\xi _{2} =\gamma _{1} (\xi _{1} )} \end{array}} \left[1-i\gamma '_{2} (x_{1} )\right]dx_{1}  -\]

\begin{equation} \label{GrindEQ__10_}
-2\int _{a_{1} }^{b_{1} }\left. \frac{\partial u(x)}{\partial x_{2} } \right|_{x_{2} =\gamma _{1} (x_{1} )} \left. \frac{\partial U(x-\xi )}{\partial x_{2} } \right|_{\begin{array}{l} {x_{2} =\gamma _{1} (x_{1} )} \\ {\xi _{2} =\gamma _{1} (\xi _{1} )} \end{array}} \left[1-i\gamma '_{1} (x_{1} )\right]dx_{1}  ,
\end{equation}

\[\left. \frac{\partial u(\xi )}{\partial \xi _{1} } \right|_{\xi _{2} =\gamma _{2} (\xi _{1} )} =\left. \frac{\partial u(\xi )}{\partial \xi _{1} } \right|_{\xi _{2} =\gamma _{1} (\xi _{1} )} -i\left. \frac{\partial u(\xi )}{\partial \xi _{2} } \right|_{\xi _{2} =\gamma _{1} (\xi _{1} )} -\]

\[-2i\int _{a_{1} }^{b_{1} }\left. \frac{\partial u(x)}{\partial x_{2} } \right|_{x_{2} =\gamma _{1} (x_{1} )} \left. \frac{\partial U(x-\xi )}{\partial x_{2} } \right|_{\begin{array}{l} {x_{2} =\gamma _{1} (x_{1} )} \\ {\xi _{2} =\gamma _{2} (\xi _{1} )} \end{array}} \left[1-i\gamma '_{1} (x_{1} )\right]dx_{1}  +\]

\begin{equation} \label{GrindEQ__11_}
+2i\int _{a_{1} }^{b_{1} }\left. \frac{\partial u(x)}{\partial x_{2} } \right|_{x_{2} =\gamma _{2} (x_{1} )} \left. \frac{\partial U(x-\xi )}{\partial x_{2} } \right|_{\begin{array}{l} {x_{2} =\gamma _{2} (x_{1} )} \\ {\xi _{2} =\gamma _{2} (\xi _{1} )} \end{array}} \left[1-i\gamma '_{2} (x_{1} )\right]dx_{1}  ,
\end{equation}

\[\left. \frac{\partial u(\xi )}{\partial \xi _{2} } \right|_{\xi _{2} =\gamma _{2} (\xi _{1} )} =2\int _{a_{1} }^{b_{1} }\left. \frac{\partial u(x)}{\partial x_{2} } \right|_{x_{2} =\gamma _{2} (x_{1} )} \left. \frac{\partial U(x-\xi )}{\partial x_{2} } \right|_{\begin{array}{l} {x_{2} =\gamma _{2} (x_{1} )} \\ {\xi _{2} =\gamma _{2} (\xi _{1} )} \end{array}} \left[1-i\gamma '_{2} (x_{1} )\right]dx_{1}  -\]

\begin{equation} \label{GrindEQ__12_}
-2\int _{a_{1} }^{b_{1} }\left. \frac{\partial u(x)}{\partial x_{2} } \right|_{x_{2} =\gamma _{1} (x_{1} )} \left. \frac{\partial U(x-\xi )}{\partial x_{2} } \right|_{\begin{array}{l} {x_{2} =\gamma _{1} (x_{1} )} \\ {\xi _{2} =\gamma _{2} (\xi _{1} )} \end{array}} \left[1-i\gamma '_{1} (x_{1} )\right]dx_{1}  .
\end{equation}
Thus, we established the following statement:

\noindent \textbf{Theorem 2.}  Let $D$  be a plane domain convex in the direction $x_{2} $, the boundary  $\Gamma $ be  Liapunov  line, then each solution of equation  \eqref{GrindEQ__1_} determined in the domain $D$ satisfies the necessary conditions \eqref{GrindEQ__8_}--\eqref{GrindEQ__12_},  containing singular integrals besides \eqref{GrindEQ__8_}.

\section*{Regularization.}

As it was said above, necessary conditions \eqref{GrindEQ__9_}--\eqref{GrindEQ__12_} contain singular addends.

\noindent  Considering   \eqref{GrindEQ__5_}, we have:

\[\left. \frac{\partial U(x-\xi )}{\partial x_{2} } \right|_{\begin{array}{l} {x_{2} =\gamma _{k} (x_{1} )} \\ {\xi _{2} =\gamma _{k} (\xi _{1} )} \end{array}} =\frac{1}{2\pi } \cdot \frac{1}{\gamma _{k} (x_{1} )-\gamma _{k} (\xi _{1} )+i(x_{1} -\xi _{1} )} =\]

\[=\frac{1}{2\pi } \frac{1}{x_{1} -\xi _{1} } \cdot \frac{1}{\gamma '_{k} (\sigma _{k} (x_{1} ,\, \xi _{1} ))+i} ,   k=1,\, 2,\]
where   $\sigma _{k} (x_{1} ,\, \xi _{1} )$  is  located  between  $x_{1} $  and  $\xi _{1} $.  Then from \eqref{GrindEQ__9_} -- \eqref{GrindEQ__12_} we find:

\[\left. \frac{\partial u(\xi )}{\partial \xi _{1} } \right|_{\xi _{2} =\gamma _{1} (\xi _{1} )} -\left. \frac{\partial u(\xi )}{\partial \xi _{1} } \right|_{\xi _{2} =\gamma _{2} (\xi _{1} )} +i\left. \frac{\partial u(\xi )}{\partial \xi _{2} } \right|_{\xi _{2} =\gamma _{2} (\xi _{1} )} =\]

\[=-\frac{1}{\pi } \int _{a_{1} }^{b_{1} }\left. \frac{\partial u(x)}{\partial x_{2} } \right|_{x_{2} =\gamma _{1} (x_{1} )} \frac{dx_{1} }{x_{1} -\xi _{1} } + ......,\]

\begin{equation} \label{GrindEQ__13_}
\left. \frac{\partial u(\xi )}{\partial \xi _{2} } \right|_{\xi _{2} =\gamma _{1} (\xi _{1} )} =\frac{i}{\pi } \int _{a_{1} }^{b_{1} }\left. \frac{\partial u(x)}{\partial x_{2} } \right|_{x_{2} =\gamma _{1} (x_{1} )} \frac{dx_{1} }{x_{1} -\xi _{1} } + ......,
\end{equation}

\[\left. \frac{\partial u(\xi )}{\partial \xi _{1} } \right|_{\xi _{2} =\gamma _{2} (\xi _{1} )} -\left. \frac{\partial u(\xi )}{\partial \xi _{1} } \right|_{\xi _{2} =\gamma _{1} (\xi _{1} )} +i\left. \frac{\partial u(\xi )}{\partial \xi _{2} } \right|_{\xi _{2} =\gamma _{1} (\xi _{1} )} =\]

\[=\frac{1}{\pi } \int _{a_{1} }^{b_{1} }\left. \frac{\partial u(x)}{\partial x_{2} } \right|_{x_{2} =\gamma _{2} (x_{1} )} \frac{dx_{1} }{x_{1} -\xi _{1} } + ......,\]

\[\left. \frac{\partial u(\xi )}{\partial \xi _{2} } \right|_{\xi _{2} =\gamma _{2} (\xi _{1} )} =-\frac{i}{\pi } \int _{a_{1} }^{b_{1} }\left. \frac{\partial u(x)}{\partial x_{2} } \right|_{x_{2} =\gamma _{2} (x_{1} )} \frac{dx_{1} }{x_{1} -\xi _{1} } + ......\]
where the sum of non-singular addends are denoted by dots.

\noindent     Considering boundary conditions \eqref{GrindEQ__2_}, from necessary conditions \eqref{GrindEQ__8_} for boundary values of the unknown function we get the following regular relations:

\[\begin{array}{l} {u(\xi _{1} ,\gamma _{1} (\xi _{1} )=u(\xi _{1} ,\, \gamma _{2} (\xi _{1} ))-} \\ {-2\int _{a_{1} }^{b_{1} }\left[\varphi _{2} (x_{1} )-\alpha _{2} u(x_{1} ,\, \gamma _{2} (x_{1} )\right]U(x_{1} -\xi _{1} ,\, \gamma _{2} (x_{1} )-\gamma _{1} (\xi _{1} ))\, \left[1-i\gamma '_{2} (x_{1} )\right]dx_{1}  +} \end{array}\]

\begin{equation} \label{GrindEQ__14_}
+2\int _{a_{1} }^{b_{1} }\left[\varphi _{1} (x_{1} )-\alpha _{1} u(x_{1} ,\, \gamma _{1} (x_{1} )\right]U(x_{1} -\xi _{1} ,\, \gamma _{1} (x_{1} )-\gamma _{1} (\xi _{1} ))\, \left[1-i\gamma '_{1} (x_{1} )\right]dx_{1}  ,
\end{equation}
In exactly the same way, from \eqref{GrindEQ__13_} we get:

\[\varphi _{1} (\xi _{1} )-\alpha _{1} u(\xi _{1} ,\, \gamma _{1} (\xi _{1} ))=\frac{i}{\pi } \int _{a_{1} }^{b_{1} }\left[\varphi _{1} (x_{1} )-\alpha _{1} u(x_{1} ,\gamma _{1} (x_{1} ))\right]\frac{dx_{1} }{x_{1} -\xi _{1} } + ......,\]

\[\varphi _{2} (\xi _{1} )-\alpha _{2} u(\xi _{1} ,\, \gamma _{2} (\xi _{1} ))=-\frac{i}{\pi } \int _{a_{1} }^{b_{1} }\left[\varphi _{2} (x_{1} )-\alpha _{2} u(x_{1} ,\gamma _{2} (x_{1} ))\right]\frac{dx_{1} }{x_{1} -\xi _{1} } + .......\]
    Finally, proceeding from \eqref{GrindEQ__14_} for boundary values of the unknown  function we get the following regular relation [9],[10]

\[\begin{array}{l} {\frac{\varphi _{1} (\xi _{1} )}{\alpha _{1} } +\frac{\varphi _{2} (\xi _{1} )}{\alpha _{2} } -\left[u(\xi _{1} ,\gamma _{1} (\xi _{1} ))+u(x_{1} ,\gamma _{2} (\xi _{1} )\right]=\frac{i}{\pi } \int _{a_{1} }^{b_{1} }\left[\frac{\varphi _{1} (x_{1} )}{\alpha _{1} } -\frac{\varphi _{2} (x_{1} )}{\alpha _{2} } \right]\frac{dx_{1} }{x_{1} -\xi _{1} }  -} \\ {-\frac{i}{\pi } \int _{a_{1} }^{b_{1} }\left\{-2\int _{a_{1} }^{b_{1} }\left[\varphi _{2} (\eta _{1} )-\alpha _{2} u(\eta _{1} ,\gamma _{2} (\eta _{1} ))\right]\,  U(\eta _{1} -x_{1} ,\, \gamma _{2} (\eta _{1} )-\gamma _{1} (x_{1} ))\left[1-i\gamma '_{2} (\eta _{1} )\right]\, d\eta _{1} +\right.  } \end{array}\]

\begin{equation} \label{GrindEQ__15_}
\left. +2\int _{a_{1} }^{b_{1} }\left[\varphi _{1} (\eta _{1} )-\alpha _{1} u(\eta _{1} ,\gamma _{1} (\eta _{1} )\right] \, U(\eta _{1} -x_{1} ,\, \gamma _{1} (\eta _{1} )-\gamma _{1} (x_{1} ))\left[1-i\gamma '_{1} (\eta _{1} )\right]\, d\eta _{1} \right\}\frac{dx_{1} }{x_{1} -\xi _{1} } +...
\end{equation}
that is regular if we interchange the integrals contained in the right hand side of \eqref{GrindEQ__15_}  and consider the singular integrals of unknown functions calculated in [11]. Thus we proved the following statement:

\textbf{Theorem 3.}  When fulfilling the conditions of Theorem 2
if $\varphi _{k} (x_{1} )$, $k=1,\, 2$ are continuously
differentiable functions vanishing at the end of the interval
$(a_{1} ,\, b_{1} )$, then \eqref{GrindEQ__15_} are regular
relations.

\section*{Fredholm property.}

Considering boundary conditions \eqref{GrindEQ__2_}, the first necessary condition \eqref{GrindEQ__8_} not containing singular integrals leads to regular relation \eqref{GrindEQ__14_}.

 Further, proceeding from boundary conditions, after regularizing two necessary conditions given in \eqref{GrindEQ__13_},  that contain singular integrals, we get a relation that has no singularity in the form \eqref{GrindEQ__15_}.

\noindent It holds :

\noindent \textbf{    Theorem 4. }When fulfilling conditions of
Theorem 3  boundary value problem
\eqref{GrindEQ__1_}--\eqref{GrindEQ__2_}  is Fredholm.

\noindent     Really, it is easy to get from \eqref{GrindEQ__14_} and \eqref{GrindEQ__15_} a system of Fredholm integral equations of second kind with respect to the unknown functions $u(x_{1} ,\gamma _{k} (x_{1} ))$, $k=1,\, 2$,  in which a kernel may have only weak singularity.

\section*{Unsolved  problems.}

\section{ The inverse problem in Tikhonov-Lavrent'ev sense. }

Let's consider the problem

\begin{equation} \label{GrindEQ__1_}
\frac{\partial ^{2} u(x)}{\partial x_{2}^{2} } +i\frac{\partial ^{2} u(x)}{\partial x_{1} \partial x_{2} } =0,\, \, \, \, x\in D,
\end{equation}

\begin{equation} \label{GrindEQ__2_}
\left. \frac{\partial u(x)}{\partial x_{2} } \right|_{x_{2} =\gamma _{k} (x_{1} )} +\alpha _{k} u(x_{1} ,\gamma _{k} (x_{1} ))=\varphi _{k} (x_{1} ),\, \, \, \, k=1,\, 2;\, \, x_{1} \in \left[a_{1} ,\, b_{1} \right],
\end{equation}
with the following complementary restriction

\[\begin{array}{l} {\left. \alpha _{1} (x_{1} )\frac{\partial u(x)}{\partial x_{2} } \right|_{x_{2} =\gamma _{1} (x_{1} )} +\left. \alpha _{2} (x_{1} )\frac{\partial u(x)}{\partial x_{2} } \right|_{x_{2} =\gamma _{2} (x_{1} )} +} \\ {+\alpha _{3} u(x_{1} ,\gamma _{1} (x_{1} ))+\alpha _{4} (x_{1} )\, u(x_{1} ,\, \gamma _{2} (x_{1} ))=\varphi _{3} (x_{1} ),\, \, \, \, \, \, x_{1} \in \left[a_{1} ,\, b_{1} \right]} \end{array}\]
where  , , $\varphi _{1} (x_{1} )$, , $k=1,\, 3,\, 4$ and  are the known,   $u(x)$, $x\in D$ è $\varphi _{2} (x_{1} )=\alpha _{2} (x_{1} )$ - are the unknown functions.

\section{ Stephan's inverse problem.}

The above mentioned boundary value problem \eqref{GrindEQ__1_},\eqref{GrindEQ__2_}, is given provided  $\alpha _{k} ,\, \, \varphi _{k} (x_{1} )$, $k=1,\, 2$, $\gamma _{1} (x_{1} )$, $\alpha _{k} (x_{1} )$, $k=\overline{1,\, 4}$ è $\varphi _{3} (x_{1} )$ - are the known,   $u(x)$, $x\in D$ è $\gamma _{2} (x_{1} )$, $x_{1} \in [a_{1} ,\, b_{1} ]$ are the unknown functions.

\section*{Reference}

 1. Courant, R. and Hilbert, D. Methods of mathematical physics, 1-2
 Interscience,
           1953 - 1962

 2. Bitsadze, A.V. Boundary value problems for second order elliptic equations, North -
            Holland, 1968

 3. Hadamard, Tohoku, Math., J., vol. 37, pp. 133-150, 1933

 4. Hadamard, L'Enseignement Math. Vol. 35, pp. 5 -42, 1936

 5. Jurayev T.D  On a boundary value problem for a composite type equation. DAN
              USSR, N4,1962,p5-8(Russian)

          6. Aliev N.A., Aliev A.M. Investigation of solutions of boundary value problems for a
              composite type equation on abounded plane domain. Dep.AzNIINTI, N688-Az from
                03.03.1987, 48 pp.

 7. Aliev N. and Jahanshahi M. Sufficient conditions for reduction of the BVP including a
               mixed PDE with non --local boundary conditions to Fredholm integral equations, INT.
               J.Math. Sci. Technol. 1997, vol. 28, N 3, 419-425

 8.  Vladimirov V.S. Equations of Mathematical Physics, Mir, 1984

9. Aliev N. and Jahanshahi M. Solation of Poisson's equation with
global, local and non -local boundary conditions, Int. J. Math.
Educ. Sci. Technol. 2002, vol.33, N2, 241-247

10. Aliev N., Hosseini S.M. An analysis of a parabolic problem
with a general (non -     local and global) supplementary linear
conditions I, II Italian Journal of Pure  and
   Applied Mathematics N12, 2002 (143-154), N13, 2003 (115-127)

 11. Gakhov F.D. Boundary value problems, Pergamon 1966

\noindent

\noindent

\end{document}